\documentclass[12pt]{amsart}
\usepackage{amsfonts}%
\usepackage{amssymb}
\usepackage{mathrsfs}
\usepackage{comment}
\usepackage[pagewise]{lineno}
\usepackage{color}
\allowdisplaybreaks

 \numberwithin{equation}{section}

\newcommand{\bsu}{\bs{u}}
\newcommand{\bso}{\bs\omega}
\newcommand{\p}{\partial}
\newcommand{\wdt}{\widetilde}
\newcommand{\ve}{\varepsilon}

\newcommand{\be}{\begin{equation}}
\newcommand{\ee}{\end{equation}}
\newcommand{\bs}{\boldsymbol}
\newcommand{\dv}{\mathrm{div}\,}
\newcommand{\cl}{\mathrm{curl}\,}
\newtheorem{Theorem}{Theorem}[section]
\newtheorem{Lemma}{Lemma}[section]
\newtheorem{Proposition}{Proposition}[section]
\newtheorem{Def}{Definition}[section]

\begin{document}

\title[Sovability of Euler equations in fractional space]
 {Solvability of Euler equations in the fractional Sobolev spaces in a bounded smooth domain}

 \author[Feng Li]{Feng Li}
\address{1. Department of Mathematics, Nanjing University, Nanjing, China}
\address{Department of Mathematics, Nagoya University, Nagoya, Japan}
 \email{lifenggoo@gmail.com}

\begin{abstract}
Euler equations are the basic system in fluid dynamics describing the motion of incompressible and inviscid ideal fluids. For a bounded smooth domain $\Omega$ in $\mathbb{R}^n$. The well-posedness of Euler equations is well-known in Sobolev spaces $W^{k,p}(\Omega)$ with the integer $k>\frac{n}{p}+1,\, 1<p<\infty$.  In this article, we study the well-posedness of Euler equations in fractional Sobolev spaces on a bounded smooth domain. We first give a priori estimates of Euler equations in fractional Hilbert-Sobolev spaces by using the energy method. For the general case of fractional Sobolev spaces, we use the characteristic method together with elliptic estimates to give similar estimates. Finally, using the a priori estimate obtained we give solvability of Euler equations in fractional Sobolev spaces. Similar to the classical case, our result is global in time in the case of two dimensions and local in the three dimensions.
\end{abstract}


\maketitle

\section{Introduction}

In this article we are concerned with the solvability of the initial-boundary value problem to the incompressible Euler equation in a bounded domain $\Omega$ in $\mathbb{R}^n, n=2,3$ with smooth boundary $\Gamma$. The Euler equations describe the motion of an incompressible inviscid fluid filling $\Omega$ as
\begin{align}
&\partial_t \bs{u}+ \bs{u} \cdot\nabla \bs{u} = -\nabla \pi \ \  \text{in}\ \Omega\times(0,T), \label{1.1}\\
&\dv\bs{u}=0 \ \  \text{in}\  \Omega\times(0,T).\label{1.2}
\end{align}
Here the unknowns are the velocity $\bs{u}(t,x)=(u_{1}, ...,u_{n})$ and the pressure $\pi=\pi(t,x)$. The system is supplemented with the initial and boundary conditions:
\begin{align}
&\bs{u}(x,0)=\bs{u}_{0}(x) \ \  \text{in} \ \Omega,\label{1.4}\\
&\bs{u}\cdot\boldsymbol{n} = 0 \ \  \text{on}\  \Gamma \times (0,T),\label{1.3}
\end{align}
where $\bs{u}_{0}=\bs{u}_{0}(x)$ is a given function and $\boldsymbol{n}$ is the unit outward normal on $\Gamma$.

The problem of well-posedness to the Euler equations has been considered for a long time by many mathematicians, starting with \cite{Lich,Wo}. In the latter article Wolibner in the first time gave the global in time existence of a regular solution in two dimensions. The notion of weak solution was first introduced by Yudovich\cite{Yu}, where he proved the global existence and uniqueness of weak solution with a bounded initial vorticity in a bounded smooth domain in $\mathbb{R}^2$. Using an abstract approach, Kato\cite{Kato1, Kato2, KP}, together with his collaborator, has proved the existence of a global solution in $\mathbb{R}^2$ and a local solution in $\mathbb{R}^3$. The existence of a local solution in the general case, i.e. $\Omega$ a domain of $\mathbb{R}^3$ with a boundary, was then proved by Ebin and Marsden in their seminal paper \cite{Ebin&Mar} by using the method of Riemannian geometry on infinite dimensional manifolds. Also, Bourguignon and Brezis\cite{Bou&Bre} have provided an alternate proof of the local existence in a more analytical way but still relying on geometric techniques. Temam\cite{Tem} applied the standard Galerkin method with a special basis to the solvability of Euler equations in a bounded smooth domain. We note that in both Bourguignon-Brezis and Temam's papers, only the solvability in the integer Sobolev spaces $W^{k,p}(\Omega)$ was considered. While in \cite{Tem1}, the main results there stated the well-posedness in $H^s(\Omega)$ with $s>\frac{n}{2}+1$ for a general real(not only integer) number $s$. However, we can not go through several key points there concerning with the fractional number $s$. Nevertheless, in the whole space or periodic domain, well-posedness in fractional Sobolev spaces or even in Besov and Triebel-Lizorkin spaces is well-known, see \cite{Lam,Chae1,Chae2,Vis1,Vis2}. We do not try to list all references here but refer to the books \cite{Chemin,Maj&Ber,MP} for detailed reports on the study of Euler equations.

The main purpose of this note is to extend the well-posedness results to the case of fractional Sobolev spaces in a bounded smooth domain. More precisely, we study the solvability of the initial-boundary value problems to Euler equations in $W^{s,p}(\Omega)$ with $1<p<\infty$ and $s\in\mathbb{R}^+, \, s>\frac{n}{p}+1$. This can be also looked as an extension of Temam's result in \cite{Tem1}. To do this we first establish a priori estimates in $W^{s,p}(\Omega)$ spaces and then use an approximation argument. The a priori estimates are based on the intrinsic characterization of fractional Sobolev spaces in a bounded domain, together with the energy method or characteristic method combing with standard estimates for elliptic equation or system. We believe these methods also work in the special case of the absence of the boundary.

The main result of this article is the well-posedness of Euler equations in the fractional Sobolev Spaces $W^{s,p}(\Omega)$, with $\Omega\subset\mathbb{R}^n$ and
\be\label{ps}
p\in (1,\infty), \ s>\frac{n}{p}+1
\ee
More precisely, we have the following theorem.
\begin{Theorem}\label{th0}
Suppose $s,p$ satisfying (\ref{ps}) and $\bs{u}_0\in W^{s,p}(\Omega)$. Then there exists a positive time $T^*$ depending on $s,p,\Omega$ and $\|\bs{u}_0\|_{s,p}$ such that there exists a unique solution $(\bs{u},\pi)$ to the initial-boundary value problem of the Euler equations. The solution belongs to the space
\[
C([0,T^*),W^{s,p}(\Omega)).
\]
Moreover, in two dimensional spaces, we have
\[
T^*=\infty,
\]
which means the solution exists for all time.
\end{Theorem}

The rest of this article is organized as follows. In the second section we give the notations used throughout the context and some preliminary results, mainly on the elliptic estimates. In section 3 we give a priori estimates for the model case $p=2$, to which the energy method can be applied. In Chapter 4 we treat the general case, where we use the method of characteristics. Finally, in the last chapter, we conclude the proof of the main results by constructing suitable approximation solutions by mollifying the initial data.


\section{Preliminaries}
In this section, we list notations and some preliminary results used throughout the rest sections.

\subsection{Basic Notations}
We denote throughout this article $\Omega$ a bounded smooth domain in $\mathbb{R}^n$, $n=2$ or $3$ with the boundary $\Gamma$. We always assume $\Gamma$ is $C^\infty$ for convenience, although the results hold for less smooth domain. The out normal vector fields on $\Gamma$ is denoted as $\bs{n}$. As the subspaces of real-valued $L^p(\Omega)$, the Sobolev spaces $W^{k,p}(\Omega)$, with $k=1,2,3,\cdots, 1\leq p \leq\infty$, are equipped with the norms
\[
\|f\|_{k,p} = \sum_{l=0}^k \|\nabla^lf\|_{p}.
\]
Hereafter we use $\nabla$ to denote the gradient operator $(\partial_{x_1}, \cdots,\partial_{x_n})$. For $p\in [1,\infty]$, $\|f\|_{p}$ is the Lebesgue $L^p(\Omega)$-norm of $f$. If $p=2$, we use $H^k(\Omega)$ to denote $W^{k,2}(\Omega)$. When $s>0$ is a real number, usually not an integer, we denote the fractional Sobolev spaces as $W^{s,p}(\Omega)$ or $H^s(\Omega)$ when $p=2$, which will be discussed in the next section.

For a function $f=f(t,x)$ defined on a time-space domain $\Omega\times [0,T]$, we use
\[
L^r(0,T,Y), \ r\in [1,\infty] \text{ and } C([0,T],Y)
\]
to denote the vector-valued space consists of functions $f$ such that $f$ is $r$-integrable and continuous in time respectively, with value in the Banach space $Y$. Usually $Y$ is a Lebesgue or Sobolev space. To distinguish a vector field from a scalar function, we use a bold symbol to denote the former.

\subsection{Fractional Sobolev Spaces}
Let $s>0$ be a real number. We denote $[s]$ the maximal integer less than or equal to $s$ and $\{s\}=s-[s]$.
The natural way to introduce the fractional Sobolev spaces is by interpolation. However, there exist two kinds of interpolation methods, real and complex interpolation, which will result in different kinds of fractional spaces. See \cite{Tri}. Nevertheless, since our method to obtain an a priori estimate(in the fractional Sobolev spaces) depends only on interpolation argument via corresponding estimates in the integer Sobolev spaces, the main results in this note hold true in any case. To make no confusion we choose the fractional Sobolev spaces via real interpolation.
\begin{Def}For $s>0$, $s\in\mathbb{N}$ and $p\in (1,\infty)$, define
\[
W^{s,p}(\Omega):= (W^{[s],p}(\Omega),W^{[s]+1,p}(\Omega))_{\{s\},p}.
\]
\end{Def}
Here of course $\Omega$ can be the whole space $\mathbb{R}^n$. In fact, in the more general setting of Besov spaces $B^s_{p,q}(\Omega)$ or Triebel-Lizorkin spaces $F^s_{p,q}(\Omega)$, $W^{s,p}(\Omega)$ coincides with $B^s_{p,p}(\Omega)$, which in turn coincides with $F^s_{p,p}(\Omega)$.

According to the extension property(see \cite{Stein}) of $W^{k,p}(\Omega)$, $k$ an integer, the fractional Sobolev spaces $W^{s,p}(\Omega)$ can be characterized via the restriction of $W^{s,p}(\mathbb{R}^n)$ to $\Omega$, that is,
\[
W^{s,p}(\Omega)=W^{s,p}(\mathbb{R}^n)|_{\Omega}.
\]


Also, it has an intrinsic characterization.  See \cite{Adams}.
\begin{Proposition}\label{prop21}
For $s>0$, the function $f\in W^{s,p}(\Omega)$ can be characterized as follows.
\[
f\in W^{[s],p}, \ [f]_{\{s\},p} := \int_{\Omega}\int_{\Omega}\frac{|\nabla^{\alpha}f(x)-\nabla^{\alpha}f(y)|^p}{|x-y|^{n+p\{s\}}} {\rm d}x{\rm d}y < \infty, \ |\alpha|=[s].
\]
\end{Proposition}
This characterization will be used only in the case of $p=2$, i.e., for $H^s(\Omega)$ spaces.

The trace space on $\Gamma$ of $W^{s,p}(\Omega)$ is $W^{s-\frac{1}{p},p}(\Gamma)$ for $s>\frac{1}{p}$, with the corresponding norm denoted as
\[
|f|_{s-\frac{1}{p},p}.
\]

According to the Sobolev's embedding theorem, for $s>\frac{n}{p}$, $W^{s,p}(\Omega)$ is embedded in $L^\infty(\Omega)$, and in fact, in $C^{s-\frac{n}{p}}(\overline{\Omega})$, the H\"{o}lderian spaces over $\Omega$. Due to this property, we have the following pointwise multiplication estimates of Moser's type.
\begin{Proposition}\label{prop22}
Suppose $s\ge 0$ and $f\in W^{s,p}(\Omega)\cap L^\infty(\Omega)$, $g\in W^{s,p}(\Omega)$. Then $fg\in W^{s,p}(\Omega)$ and there exists a constant $C>0$ independent of $f,g$ such that
\[
\|fg\|_{s,p} \leq \left( \|f\|_{\infty} + \|f\|_{s,p}\right)\|g\|_{s,p}.
\]
Moreover, if $s>\frac{n}{p}$ and $f,g\in W^{s,p}(\Omega)$, then $fg\in W^{s,p}(\Omega)$ and there exists a constant $C$ independent of $f,g$ such that
\[
\|fg\|_{s,p} \leq C\left(\|f\|_{\infty}\|g\|_{s,p} + \|f\|_{s,p}\|g\|_{\infty}\right) \leq C\|f\|_{s,p}\|g\|_{s,p}.
\]
\end{Proposition}
The following proposition on composition of functions is a consequence of this multiplication property.

\begin{Proposition}\label{prop23}
Let $f\in W^{s,p}(\Omega), s>1+\frac{n}{p}$, $f(0)=0$ and $X:\Omega\rightarrow\Omega$ a one-to-one mapping in $\Omega$ with every component $\nabla X_i \in W^{s-1,p}(\Omega),i=1,\cdots,n$. Then $f(X)\in W^{s,p}(\Omega)$ and
\[
\|f(X)\|_{s,p} \leq C\|f\|_{s,p}\|\nabla X\|_{s-1,p}.
\]
\end{Proposition}

At the end of this section we give an embedding result of Bre\'{z}is-Wainger(\cite{BW}). For the present form one can refer to \cite{SWZ}. First we need to introduce John-Nirenberg's space of bounded mean oscillation, the $BMO(\Omega)$ space. To this end let us denote the mean value of a function $f$ on a measurable set $D$ as
\[
f_{D}=\frac{1}{|D|}\int_D f(y){\rm d}y,
\]
where $|D|$ is the Lebesgue measure of $D$. A ball in $\mathbb{R}^n$ is denoted as $B$. The $BMO(\Omega)$ space is a subspace of $L^2(\Omega)$ consists of functions $f$ such that
\[
[f]_* := \sup_{B} \frac{1}{|B\cap\Omega|}\int_{B\cap\Omega} |f(y)-f_{B\cap\Omega}|{\rm d}y< \infty,
\]
equipped with the norm
\[
\|f\|_* = \|f\|_2 + [f]_*.
\]
Obviously $L^{\infty}(\Omega)\subset BMO(\Omega)$ and it is well known that $W^{\frac{n}{p},p}(\Omega)$ is continuously embedded in $BMO(\Omega)$. See, for example, \cite{Tay1}. The following proposition gives a log-type embedding of $W^{s,p}(\Omega)$ with $s>\frac{n}{p}$ into $L^{\infty}(\Omega)$ with the help of $BMO$-norm.
\begin{Proposition}\label{prop235}
Suppose
\[
s> \frac{n}{p}.
\]
There exists a constant $C$ depending only on $n,s,p$ and $\Omega$ such that for any $f\in W^{s,p}(\Omega)$,
\[
\|f\|_{\infty} \leq C \|f\|_{*}\ln(1+\|f\|_{s,p}).
\]
\end{Proposition}
Note that this log-type estimate hold true for any dimension number $n$ but we only need it in $2D$ case. The proof of this proposition is based on the embedding of $W^{s,p}(\Omega)$ into $C^{s-\frac{n}{p}}(\overline{\Omega})$, together with the following estimate for $BMO$ functions.
\[
|f_{B_r\cap\Omega}| \leq C(1+|\ln{r}|)\|f\|_*,
\]
where $B_r$ is a ball with radius $r\leq 1$. We refer to \cite{SWZ} for more details.

\subsection{Estimates for Elliptic Equations}

In this section, we list some classical and well-known estimates for elliptic equations. We do not give exact references here but refer to \cite{Gil&Tru} as well as \cite{AgmDou&Nir} for modern and general treatment. Usually, in these classical works, only the integer Sobolev spaces are considered. But one can get the corresponding estimates in fractional Sobolev spaces via interpolation argument. In fact, the elliptic estimates in the most general case are well developed by Triebel in \cite{Tri} in the most general setting, that is, in spaces of Besov($B^s_{p,q}(\Omega)$) and Triebel-Lizorkin($F^s_{p,q}(\Omega)$) type, almost without any restriction on the parameters $p,q\in (0,\infty], \ s\in \mathbb{R}$.

Although we only use them in dimensional two and three, these estimates given below are in fact true in any dimensions. We begin with the Dirichlet problems for the Laplace equation. Consider
\be\label{el1}
\Delta{v} = f \ \text{in} \ \Omega, \ \ v=g \text{ on }\Gamma.
\ee
\begin{Proposition}\label{prop24}
Suppose $v$ is the solution of (\ref{el1}) with $f\in W^{s,p}(\Omega),$
$ g\in W^{s+2-\frac{1}{p},p}(\Gamma)$ with $s\ge -1,p\in (1,\infty)$.
Then there exists a constant depending only on $\Omega$ and $n,s,p$ such that
\[
\|v\|_{s+2,p} \leq C \left(\|f\|_{s,p} + |g|_{s+2-\frac{1}{p},p} \right).
\]
\end{Proposition}
Similarly for the Neumann problem
\be\label{el2}
\Delta{v} = f \ \text{in} \ \Omega, \ \ \frac{\partial v}{\partial \bs{n}}=g \text{ on }\Gamma,
\ee
we have
\begin{Proposition}\label{prop25}
Suppose $v$ is the solution of (\ref{el2}) with $f\in W^{s,p}(\Omega), $
$g\in W^{s+1-\frac{1}{p},p}(\Gamma)$ with $s\ge 0,p\in (1,\infty)$.
Then there exists a constant depending only on $\Omega$ and $n,s,p$ such that
\[
\|\nabla v\|_{s+1,p} \leq C \left(\|f\|_{s,p} + |g|_{s+1-\frac{1}{p},p} \right).
\]
\end{Proposition}
Next we consider the following "div-curl" system.
\be\label{el3}
\dv\bs{v} = 0, \ \cl\bs{v} = \bs{f} \text{ in }\Omega, \ \bs{v}\cdot\bs{n} = 0 \text{ on }\Gamma.
\ee
\begin{Proposition}\label{prop27}
Suppose $\bs{v}$ is the solution of (\ref{el3}) with $\bs{f}\in W^{s,p}(\Omega)$ with $s\ge 0,p\in (1,\infty)$. Then there exists a constant depending only on $\Omega$ and $s,p$ such that
\[
\|\bs{v}\|_{s+1,p} \leq C \|\bs{f}\|_{s,p}.
\]
\end{Proposition}
In case of $s=k$ an integer, see \cite{Wahl} for a proof. The general case follows from interpolation.

In treating $2D$ Euler equations we need an $BMO$ estimate for Dirichlet problem. We formulate it as follows.
\begin{Proposition}\label{prop28}
Let $v$ be the solution of the Dirichlet problem (\ref{el1}) with $g=0$ and $f\in BMO(\Omega)$. Then
$\nabla^2v\in BMO(\Omega)$ and there is a constant $C$ independent of $f$ such that
\[
\|\nabla^2v\|_{*} \leq C\|f\|_{*}.
\]
\end{Proposition}
See \cite{Aq} for proof in more general settings.

\subsection{Gr\"{o}nwall's Inequality and Transport Equations}
In this section we give estimates for transport equations in the fractional Sobolev spaces via interpolation argument. To the author's knowledge, this kind of estimates has been used in \cite{KP} with the underlying space as $\mathbb{R}^n$. We begin with the well-known inequality of Gr\"{o}nwall's type.
\begin{Lemma}\label{le21}
Let $g(t)$ be a nonnegative continuous function defined on $[0,T]$. Suppose there exists two nonnegative functions $h(t), c(t)\in L^1(0,T)$ such that
\[
g'(t) \leq h(t) + c(t)g(t).
\]
Then
\[
g(t) \leq  \left( g(0) + \int_0^t h(\tau){\rm d}\tau \right)\exp\left(\int_0^t c(\tau){\rm d}\tau\right)
\]
for any $t\in [0,T]$.
\end{Lemma}
The following simple lemma concerns the case of nonlinear increasing of $g'$.
\begin{Lemma}\label{le22}
Let $g(t)$ be a nonnegative continuous function defined on $[0,T]$. Suppose there exists two constant $C>0$ and $m >1$ such that
\[
g'(t) \leq Cg^m(t).
\]
Then there is a $T^*\leq T$ depending on $C, m$ and $g(0)$ such that for $t\in [0,T^*]$,
\[
g(t) \leq 2g(0).
\]
\end{Lemma}

Now consider the following linear transport equation in $\Omega$.
\be\label{trans}
\partial_t u + \bs{v}\cdot\nabla{u} = f, \ u(x,0)=u_0(x).
\ee
Here $f=f(t,x)$ and $\bs{v}=\bs{v}(t,x)$ are known function and vector field defined in $\Omega\times[0,T]$ respectively. Moreover, $\bs{v}$ satisfies
\[
\dv\bs{v} = 0, \  \bs{v}\cdot\bs{n} = 0 \text{ on }\Gamma.
\]
for all $t\in [0,T]$.
\begin{Proposition}\label{prop29}
Given $u_0\in W^{s,p}(\Omega)$, $s\in [0,1], \ p\in (1,\infty)$. Now we suppose $\bs{v}\in L^1(0,T; W^{1,\infty}(\Omega))$ and $f\in L^1(0,T,W^{s,p}(\Omega))$. Let $u(t,x)\in C([0,T],W^{s,p}(\Omega))$ be the solution to (\ref{trans}). Then there exists a constant $C$ depending only on $s,p$ such that for any $t\in [0,T]$,
\be\label{tr1}
\|u(t,\cdot)\|_{s,p} \leq C\left( \|u_0\|_{s,p}+ \int_0^t \|f(\tau,\cdot)\|_{s,p}{\rm d}\tau \right)\exp\left(C\int_0^\tau\|\nabla\bs{v}(\tau,\cdot)\|_{\infty}{\rm d}\tau\right).
\ee
In the case of $p=\infty$ and $f\equiv 0$ we have the following obvious estimate.
\be\label{tr2}
\|u(t,\cdot)\|_{\infty} \leq \|u_0\|_{\infty}.
\ee
\end{Proposition}
\noindent{\bf Proof:} The proof is standard. Since the equation is linear, we only need to prove the case of $s=0,1$. Estimate (\ref{tr1}) then follows by interpolation argument. To this end we note that by the condition imposed on $\bs{v}$, it follows that for any $p\in (1,\infty)$,
\[
\frac{{\rm d}}{{\rm d}t} \|u(t,\cdot)\|^p_{p} \leq p\|f(t,\cdot)\|_{p}\|u(t,\cdot)\|^{p-1}_{p}.
\]
Consequently,
\[
\frac{{\rm d}}{{\rm d}t} \|u(t,\cdot)\|_{p} \leq p\|f(t,\cdot)\|_{p}.
\]
Thus for any $t\in [0,T]$,
\be\label{tr3}
\|u(t,\cdot)\|_{p} \leq \|u_0\|_{p} + p\int_0^t \|f(\tau,\cdot)\|_{p}{\rm d}\tau
\ee
If $f\equiv 0$, then
\[
\|u(t,\cdot)\|_{p} \leq \|u_0\|_{p}.
\]
Estimate (\ref{tr2}) follows by taking $p\rightarrow\infty$. Now taking derivatives with respect to $x$ yields
\[
\partial_t\nabla{u} + \bs{v}\cdot\nabla\nabla{u}= \nabla{f}(t,x) - \nabla\bs{v}\cdot\nabla{u}.
\]
Hence we have
\[
\frac{{\rm d}}{{\rm d}t} \|\nabla{u}(t,\cdot)\|^p_{p} \leq p\|f(t,\cdot)\|_{p}\|\nabla{u}(t,\cdot)\|^{p-1}_{p} + p\|\nabla\bs{v}\|_{\infty}\|\nabla{u}(t,\cdot)\|^p_{p}.
\]
Applying Gr\"{o}nwall's inequality gives the desired estimate (\ref{tr1}) for $s=1$.

\section{A Priori Estimates in $H^s(\Omega)$}
In this section, the energy method is applied to obtain a priori estimates for Euler equations in the fractional Hilber-Sobolev spaces $H^s(\Omega)$. The intrinsic characterization of $H^s(\Omega)$ makes it possible to apply this method. It will be seen that, compared to the characteristic method(combing with energy method), the use of energy method is rather complex in our setting. We discuss the $2D$ and $3D$ cases separately. In this section we use $\|f\|_{s}$ to denote the $H^s(\Omega)$-norm for notation brevity.

\subsection{A Priori Estimates in $2D$}
In $2D$ space, the vorticity $\omega$ is a scalar function, which is transported by the corresponding velocity.
\be\label{e31}
\partial_t\omega + \bs{u}\cdot\nabla\omega = 0, \ \omega(x,0) = \omega_0(x) \text{ in }\Omega.
\ee
The velocity $\bs{u}$ satisfies
\be\label{e311}
\dv\bs{u} = 0, \ \cl\bs{u}=\omega, \ \bs{u}(x,0)=\bs{u}_0(x), \ \text{ in }\Omega, \ \bs{u}\cdot\bs{n} = 0 \text{ on }\Gamma.
\ee
Moreover, by introducing the stream function $\psi(t,x)$ such that
\[
\bs{u}(t,x)= \nabla^\bot\psi := (\partial_{x_2}\psi(t,x), -\partial_{x_1}\psi(t,x)),
\]
we have
\[
-\Delta\psi = \omega, \, \text{in}\ \Omega, \, \psi(t,x)=0 \ \text{on} \ \Gamma.
\]
By Proposition \ref{prop24} and \ref{prop28},
\[
\|\psi(t,\cdot)\|_{s+1,p} \leq C \|\omega(t,\cdot)\|_{s-1,p}, \ \|\nabla^2\psi(t,\cdot)\|_{*}\leq C\|\omega(t,\cdot)\|_{*}\leq C\|\omega(t,\cdot)\|_{\infty}
\]
with some constant $C$ independent of $t$ and $\omega$. Consequently,
\be\label{e30}
\|\bs{u}(t,\cdot)\|_{s,p} \leq C \|\omega(t,\cdot)\|_{s-1,p}, \ \|\nabla\bs{u}(t,\cdot)\|_{*}\leq C\|\omega(t,\cdot)\|_{\infty}.
\ee
With these preparation we are ready to prove the following a priori estimates for $\omega$ in $H^{s}(\Omega)$ space.
\begin{Proposition}\label{prop31}
Assume $\bs{u}_0\in H^{s+1}(\Omega)$ with $1<s<2$ and $\bs{u}(t,x),\omega(t,x)$ is the solution to (\ref{e31})-(\ref{e311}). There exists a positive constant $C$ depending only on $s,\Omega$ and $\bs{u}_0$ such that for any $t>0$,
\be\label{e34}
\|\omega(t,\cdot)\|_{s} \leq C\exp(\exp(Ct)), \ \|\bs{u}(t,\cdot)\|_{s+1} \leq C\exp(\exp(Ct)).
\ee
\end{Proposition}
\noindent \textbf{Proof}: First note that we have
\[
\|\omega(t,\cdot)\|_{L^{2}(\Omega)} \leq \|\omega_0\|_{L^{2}(\Omega)} \leq C\|\bs{u}_0\|_{1}
\]
as well as
\be\label{e35}
\|\omega(t,\cdot)\|_{L^{\infty}(\Omega)} \leq \|\omega_0\|_{L^{\infty}(\Omega)} \leq C\|\bs{u}_0\|_{s}
\ee
due to Sobolev's embedding lemma.

To estimate $\|\omega\|_{s}$-norm by using the intrinsic characterization of $H^s(\Omega)$ we first write (\ref{e31}) in two independent variables $x$ and $y$.
\[
\begin{split}
    &\partial_t{\omega}(t,x)+\bs{u}(t,x)\cdot\nabla_x{\omega}(t,x) = 0,\ x\in\Omega,\\
    &\partial_t{\omega}(t,y)+\bs{u}(t,y)\cdot\nabla_y{\omega}(t,y) = 0,\ x\in\Omega.
\end{split}
\]

It follows that
\[
\begin{split}
    &\partial_t(\omega(t,x)-\omega(t,y))+\bs{u}(t,x)\cdot(\nabla_x\omega(t,x)-\nabla_y\omega(t,y))\\
    =&-(\bs{u}(t,x)-\bs{u}(t,y))\cdot\nabla_y\omega(t,y).
\end{split}
\]

Let us denote $\delta=s-1\in(0,1)$ temporarily. Multiplying the equation by
$$
\frac{\omega(t,x)-\omega(t,y)}{|x-y|^{2+2\delta}}
$$
and integrating on $\Omega\times\Omega$ give rise to
\[
\begin{split}
    &\frac{1}{2}\frac{{\rm{d}}}{{\rm{d}}t}\int_{\Omega\times\Omega} \frac{|\omega(t,x)-\omega(t,y)|^2}{|x-y|^{2+2\delta}}{\rm{d}} x{\rm{d}} y\\
    = & -\frac{1}{2}\int_{\Omega\times\Omega}(\bs{u}(t,x)\cdot\nabla_x|x-y|^{-2-2\delta})
|\omega(t,x)-\omega(t,y)|^2{\rm{d}}{x}{\rm{d}}{y}.
\end{split}
\]

Switching the position of $x$ and $y$, then adding it to the former one gives rise to
\[
\begin{split}
    &\frac{{\rm{d}}}{{\rm{d}}t}\int_{\Omega\times\Omega}\frac{|\omega(t,x)-\omega(t,y)|^2}{|x-y|^{2+2\delta}}
{\rm{d}}{x}{\rm{d}}{y}\\
=&\frac{1}{2}\int_{\Omega\times\Omega}\left(\bs{u}(t,x)\cdot\nabla_x\frac{1}{|x-y|^{2+2\delta}}
+\bs{u}(t,y)\cdot\nabla_y\frac{1}{|x-y|^{2+2\delta}}\right)\\
&\qquad\times|\omega(t,x)-\omega(t,y)|^2{\rm{d}}{x}{\rm{d}}{y}.
\end{split}
\]

The term on the right-hand side is controlled by
\[
\begin{split}
    & C\int_{\Omega\times\Omega}\left|\bs{u}(t,x)-\bs{u}(t,y)\right|
\frac{|\omega(t,x)-\omega(t,y)|^2}{|x-y|^{3+2\delta}}{\rm{d}}{x}{\rm{d}}{y}\\
= & C \int_{\Omega\times\Omega}\frac{\left|\bs{u}(t,x)-\bs{u}(t,y)\right|}{|x-y|}
\frac{|\omega(t,x)-\omega(t,y)|^2}{|x-y|^{2+2\delta}}{\rm{d}}{x}{\rm{d}}{y}\\
\leq & C \|\nabla\bs{u}(t)\|_{L^{\infty}}\int_{\Omega\times\Omega}
\frac{|\omega(t,x)-\omega(t,y)|^2}{|x-y|^{2+2\delta}}{\rm{d}}x{\rm{d}}y.
\end{split}
\]

It follows that
$$
\frac{{\rm{d}}}{{\rm{d}}t}\|\omega(t,\cdot)\|_{s-1} \leq C\|\nabla\bs{u}(t,\cdot)\|_{L^\infty}\|\omega(t,\cdot)\|_{s-1}.
$$
By Gr\"{o}nwall's inequality,
\be\label{e36}
\|\omega(t,\cdot)\|_{s-1} \leq C\exp{\int_0^t\|\nabla\bs{u}(t,\cdot)\|_{L^\infty}}{\rm{d}}\tau,
\ee
for any $t>0$.

To continue we differentiate both sides of (\ref{e31}) with respect to $x=x_1,x_2$ to give
$$
\partial_t\partial_x\omega(t,x)+\bs{u}(t,x)\cdot\nabla\partial_x\omega(t,x) = -\partial_x\bs{u}(t,x)\cdot\nabla\omega(t,x):=f(t,x).
$$
Similarly,
$$
\partial_t\partial_y\omega(t,y)+\bs{u}(t,y)\cdot\nabla\partial_y\omega(t,y) = -\partial_y\bs{u}(t,y)\cdot\nabla\omega(t,y):=f(t,y).
$$
Now repeat the procedure as before, noting $\dv\bs{u}=0, \bs{u}\cdot\bs{n}=0$, and after integration by parts we arrive at the following energy equality for $\partial\omega$.
\[
\begin{split}
    &\frac{{\rm{d}}}{{\rm{d}}t}\int_{\Omega\times\Omega}
\frac{|\partial_x\omega(t,x)-\partial_y\omega(t,y)|^2}{|x-y|^{2+2\delta}}{\rm{d}}{x}{\rm{d}}{y}\\
= & \int_{\Omega\times\Omega}
(\bs{u}(t,x)\cdot\nabla_x|x-y|^{-2-2\delta}+\bs{u}(t,y)\cdot\nabla_y|x-y|^{-2-2\delta})\\
&\qquad\times|\partial_x\omega(t,x)-\partial_y\omega(t,y)|^2{\rm{d}}{x}{\rm{d}}{y}\\
&+2\int_{\Omega\times\Omega}\left(f(t,x)-f(t,y)\right)
\frac{\partial_x\omega(t,x)-\partial_y\omega(t,y)}{|x-y|^{2+2\delta}}{\rm{d}}{x}{\rm{d}}{y}.
\end{split}
\]

The first term on the right-hand side is controlled by(recall that $\delta=s-1$)
$$
\|\nabla\bs{u}\|_{L^{\infty}(\Omega)}\|\nabla\omega\|^2_{s-1}
$$
by a similar argument as before. The second term is controlled by
$$
\|f(t,\cdot)\|_{s-1}\|\nabla\omega(t,\cdot)\|_{s-1}.
$$
To control $\|f\|_{s-1}$ we need divergence free condition of $\bs{u}$. Note that for $\partial_x, x=x_1,x_2$,
$$
f(t,x)=\partial_x\bs{u}\cdot\nabla\omega = \dv(\omega\partial_x\bs{u}).
$$
Hence by Proposition \ref{prop22} and \ref{prop24},
\[
\begin{split}
   & \|f(t,\cdot)\|_{s-1} \leq \|\omega\partial_x\bs{u}\|_{s} \leq C \left( \|\omega\|_{L^\infty}\|\partial_x\bs{u}\|_{s} + \|\partial_x\bs{u}\|_{L^\infty}\|\omega\|_{s}\right)\\
   \leq & C \|\nabla\bs{u}\|_{L^\infty} \|\omega\|_{s}.
\end{split}
\]

Collecting all these estimates we finally arrived at
$$
\frac{{\rm{d}}}{{\rm{d}}t}\|\omega(t,\cdot)\|_{s} \leq C\|\nabla\bs{u}(t,\cdot)\|_{L^\infty}\|\omega(t,\cdot)\|_{s}.
$$
By Gr\"{o}nwall's inequality,
\be\label{e37}
\|\omega(t,\cdot)\|_{s} \leq C\exp{\int_0^t\|\nabla\bs{u}(t,\cdot)\|_{L^\infty}}{\rm{d}}\tau,
\ee
for any $t>0$.

To close estimate we use (\ref{e30}) and (\ref{e35}) to find
\[
\begin{split}
    &\|\nabla\bs{u}(t,\cdot)\|_{L^\infty} \leq C\|\nabla\bs{u}(t,\cdot)\|_* \ln(1+\|\nabla\bs{u}(t,\cdot)\|_{s})\\
    \leq & C\|\omega(t,\cdot)\|_{\infty}\ln(1+\|\omega(t,\cdot)\|_{s}) \leq C \ln(1+\|\omega(t,\cdot)\|_{s}),
\end{split}
\]
where in the first inequality we used Proposition \ref{prop235}.
Substituting this into (\ref{e36}) yields
$$
\|\omega(t,\cdot)\|_{s} \leq C\exp{\left(C\int_0^t \ln\left( 1 + \|\omega(\tau,\cdot)\|_{s} \right) {\rm{d}}\tau \right)}.
$$
Hence,
$$
\ln \left( 1+ \|\omega(t,\cdot)\|_{s}\right) \leq  C \int_0^t \ln \left( 1+ \|\omega(t,\cdot)\|_{s}\right)
$$
The proof of (\ref{e34}) for $\omega$ is concluded by Gr\"{o}nwall's inequality while the corresponding estimates for $\bs{u}, \pi$ follow from elliptic estimates.

\subsection{A Priori Estimates in $3D$}
In $3D$, the evolution of vorticity $\bs\omega(t,x)$ is formulated as
\be\label{e300}
\partial_t\bs{\omega} + \bs{u}\cdot\nabla\bs{{\omega}} = \bs\omega\cdot\nabla\bs{u}, \, \bs\omega(x,0)=\bs\omega_0(x).
\ee
While the velocity satisfy
\be\label{e3000}
\dv\bs{u} = 0, \ \cl\bs{u} = \bs\omega \text{ in }\Omega, \ \bs{u}\cdot\bs{n} = 0 \text{ on }\Gamma.
\ee
By Proposition \ref{prop27}, for any $s\ge 1$,
\be\label{e301}
\|\bs{u}(t,\cdot)\|_{s} \leq C \|\bs\omega(t,\cdot)\|_{s-1}
\ee
with the positive constant $C$ independent of $t$ and $\bs\omega$. We have the following local estimates for $3D$ Euler equations in the fractional Sobolev spaces $H^s(\Omega)$.
\begin{Proposition}\label{prop32}
Assume $s\in (\frac{5}{2},3)$, $\boldsymbol{\bs{u}_0}\in H^{s}(\Omega)$. There exists a positive time $T>0$ and a constant $C>0$ depending on $\|\bs{u}_0\|_{s}$ such that for any $t\in [0,T]$,
\be\label{e37}
\|\bs{\omega}(t,\cdot)\|_{s-1}\leq C, \, \|\bs{u}(t,\cdot)\|_{s} + \|\pi(t,\cdot)\|_{s} \leq C.
\ee
\end{Proposition}

\noindent \textbf{Proof}. Similar as before we consider the difference between $\bs\omega(t,x)$ and $\bs\omega(t,y)$.
\[
\begin{split}
    &\partial_t(\bs{\omega}(t,x)-\bs{\omega}(t,y))+\bs{u}(t,x)\cdot(\nabla_x\bs{\omega}(t,x)-
\nabla_y\bs{\omega}(t,y))\\
 =\ & \bs{\omega}(t,x)\cdot\nabla_x\bs{u}(t,x)-\bs{\omega}(t,y)\cdot\nabla_y\bs{u}(t,y)
-(\bs{u}(t,x)-\bs{u}(t,y))\cdot\nabla_y\bs{\omega}(t,y).
\end{split}
\]

Let $\delta\in (0,1)$. Multiplying the equation by $\frac{\bs{\omega}(t,x)-\bs{\omega}(t,y)}{|x-y|^{3+2\delta}}$ and integrating on $\Omega\times\Omega$ gives rise to
\[
\begin{split}
    &\frac{1}{2}\frac{{\rm{d}}}{\rm{d}t}\int_\Omega \int_\Omega\frac{|\bs{\omega}(t,x)-\bs{\omega}(t,y)|^2}{|x-y|^{3+2\delta}}{\rm{d}}{x}{\rm{d}}{y}\\
    = & \int_{\Omega}\int_{\Omega}\dv_x\bs{u}(t,x)
\frac{|\bs{\omega}(t,x)-\bs{\omega}(t,y)|^2}{|x-y|^{3+2\delta}}{\rm{d}}{x}{\rm{d}}{y}\\
&\qquad -\frac{1}{2}\int_{\Omega}\int_{\Omega}(\bs{u}(t,x)\cdot\nabla_x\frac{1}{|x-y|^{3+2\delta}})
|\bs{\omega}(t,x)-\bs{\omega}(t,y)|^2{\rm{d}}{x}{\rm{d}}{y}\\
&\qquad +\int_{\Omega}\int_{\Omega}\left(\bs{\omega}(t,x)\cdot\nabla\bs{u}(t,x)-
\bs{\omega}(t,y)\cdot\nabla\bs{u}(t,y)\right)\cdot
\frac{\bs{\omega}(t,x)-\bs{\omega}(t,y)}{|x-y|^{3+2\delta}}{\rm{d}}{x}{\rm{d}}{y}.
\end{split}
\]

Using $\dv\bs{u}=0$,
\[
\begin{split}
    &\frac{1}{2}\frac{{\rm{d}}}{\rm{d}t} \int_\Omega \int_\Omega\frac{|\bs{\omega}(t,x)-\bs{\omega}(t,y)|^2}{|x-y|^{3+2\delta}}{\rm{d}}{x}{\rm{d}}{y}\\
    = & -\frac{1}{2}\int_{\Omega}\int_{\Omega}\left(\bs{u}(t,x)\cdot\nabla_x\frac{1}{|x-y|^{3+2\delta}}\right)
|\bs{\omega}(t,x)-\bs{\omega}(t,y)|^2{\rm{d}}{x}{\rm{d}}{y}\\
&+\int_{\Omega}\int_{\Omega}\left(\bs{\omega}(t,x)\cdot\nabla\bs{u}(t,x)-
\bs{\omega}(t,y)\cdot\nabla\bs{u}(t,y)\right)\cdot
\frac{\bs{\omega}(t,x)-\bs{\omega}(t,y)}{|x-y|^{3+2\delta}}{\rm{d}}{x}{\rm{d}}{y}.
\end{split}
\]

It follows that
\begin{equation}\label{e38}
\begin{split}
    &\frac{{\rm{d}}}{{\rm{d}}t}\int_\Omega \int_\Omega\frac{|\bs{\omega}(t,x)-\bs{\omega}(t,y)|^2}{|x-y|^{3+2\delta}}{\rm{d}}{x}{\rm{d}}{y}\\
    = & - \int_{\Omega}\int_{\Omega}\big{(}\bs{u}(t,x)\cdot\nabla_x\frac{1}{|x-y|^{3+2\delta}}
+\bs{u}(t,y)\cdot\nabla_y\frac{1}{|x-y|^{3+2\delta}}\big{)}\\
&\qquad\times|\bs{\omega}(t,x)-\bs{\omega}(t,y)|^2{\rm{d}}{x}{\rm{d}}{y}\\
& + 2\int_{\Omega}\int_{\Omega}(\bs{\omega}(t,x)\cdot\nabla\bs{u}(t,x)
-\bs{\omega}(t,y)\cdot\nabla\bs{u}(t,y))\cdot
\frac{\bs{\omega}(t,x)-\bs{\omega}(t,y)}{|x-y|^{3+2\delta}}{\rm{d}}{x}{\rm{d}}{y}.
\end{split}
\end{equation}

The first term on the right-hand side is controlled by
\[
\begin{split}
    &\int_{\Omega}\int_{\Omega}\frac{\left|\bs{u}(t,x)-\bs{u}(t,y)\right|} {|x-y|^{4+2\delta}}|\bs{\omega}(t,x)-\bs{\omega}(t,y)|^2{\rm{d}}{x}{\rm{d}}{y}\\
    \leq & C \|\nabla\bs{u}(t,\cdot)\|_{L^{\infty}}\int_{\Omega}\int_{\Omega}
\frac{|\bs{\omega}(t,x)-\bs{\omega}(t,y)|^2}{|x-y|^{3+2\delta}}{\rm{d}}{x}{\rm{d}}{y}\\
\leq & C\|\nabla\bs{u}(t,\cdot)\|_{L^\infty}\|\bs{\omega}(t,\cdot)\|^2_{\delta}.
\end{split}
\]

The second term on the right-hand side of (\ref{e38}) is controlled by
\[
\begin{split}
    &C\|(\bs\omega\cdot\nabla\bs{u})(t,\cdot)\|_{\delta} \|\bs\omega(t,\cdot)\|_{\delta}\\
    \leq & C\left(\|(\bs\omega(t,\cdot)\|_{L^\infty}\|\nabla\bs{u}(t,\cdot)\|_{\delta} + \|\bs\omega(t,\cdot)\|_{\delta} \|\nabla\bs{u}(t,\cdot)\|_{L^\infty} \right) \|\bs\omega(t,\cdot)\|_{\delta}\\
    \leq & C \|\nabla\bs{u}(t,\cdot)\|_{L^\infty} \|\bs\omega(t,\cdot)\|^2_{\delta},
\end{split}
\]
according to multiplication property and elliptic estimates (\ref{e301}). We thus conclude that
\be\label{e39}
\frac{{\rm d}}{{\rm d}t} \|\bs{\omega}(t,\cdot)\|_{s-2} \leq C \|\nabla\bs{u}(t,\cdot)\|_{L^\infty}\|\bs{\omega}(t,\cdot)\|_{s-2}
\ee
by taking $\delta=s-2 \in ({1\over 2}, 1)$.

We need further to give estimates on $\|\nabla\omega\|_{s-2}$.
By taking spatial derivative $\partial_x,x=x_1,x_2,x_3$ on both side of (\ref{e300}),
\[
\partial_t\partial_x\bso + \bsu\cdot\p_x\nabla\bso = - \p_x\bsu\cdot\nabla\bso + \p_x(\bs\omega(t,x)\cdot\nabla\bs{u}(t,x)):= \bs{f}(t,x).
\]
Note that since
\[
\p_x\bsu\cdot\nabla\bso = \dv(\p_x\bsu \otimes \bso ),
\]
\[
\bs{f}=  \dv(\p_x\bsu \otimes \bso ) + \p_x \left( \bs\omega\cdot\nabla\bs{u}\right).
\]
Consequently, by $s-1>\frac{3}{2}$ and the multiplication property,
\[
\begin{split}
    &\|\bs{f}(t,\cdot)\|_{s-2} \leq \|(\p_x\bsu \otimes \bso)(t,\cdot)\|_{s-1}+ \| (\bs\omega\cdot\nabla\bs{u})(t,\cdot) \|_{s-1}\\
    \leq & C \left( \|\p_x\bsu\|_{L^\infty} \|\bso\|_{s-1} + \|\p_x\bsu\|_{s-1}\| \bso \|_{L^\infty} \right)\\
    &\qquad+\left( \|\bso\|_{s-1}\|\nabla\bsu\|_{L^\infty} + \|\bso\|_{L^\infty}\|\nabla\bsu\|_{s-1} \right)\\
    \leq & C \|\nabla\bsu(t,\cdot)\|_{L^\infty} \|\bso(t,\cdot)\|_{s-1}.
\end{split}
\]

With this estimate for $\bs{f}$ and exactly the same argument as before we have
\be\label{e310}
\frac{{\rm d}}{{\rm d}t}\|\nabla\bso(t,\cdot)\|_{s-2} \leq C\|\nabla\bsu(t,\cdot)\|_{L^\infty} \|\bso(t,\cdot)\|_{s-1}
\ee
Combining with (\ref{e39}) yields
\be\label{e312}
\frac{{\rm d}}{{\rm d}t}\|\bso(t,\cdot)\|_{s-1} \leq C\|\nabla\bsu(t,\cdot)\|_{L^\infty} \|\bso(t,\cdot)\|_{s-1},
\ee
together with the embedding
\[
\|\nabla\bsu\|_{L^\infty} \leq C \|\bsu\|_{s} \leq C \|\bso\|_{s-1},
\]
we conclude the proof of (\ref{e37}) for $\|\bso(t,\cdot)\|_{s-1}$ by Lemma \ref{le22}. Estimates for $\bsu$ and $\pi$ follow from elliptic estimates.

\noindent{\bf Remark}: In Proposition \ref{prop31} and \ref{prop32}, we only give a priori estimates in $H^s(\Omega)$ with $s<3$. For larger $s$, similar estimates can be obtained with the help of multiplication properties and elliptic estimates.

\section{A Priori Estimates in $W^{s,p}(\Omega)$}
In this section, we use the method of characteristics to deduce a priori estimates of the solution $(\bs{u},\pi)$ to Euler equations in the fractional Sobolev spaces $W^{s,p}(\Omega)$, with $p \in (1,\infty), s>\frac{n}{p} + 1$. Neither energy method nor Littlewood-Paley analysis can be applied, due to the facts that $p\neq 2$ and $\Omega\neq \mathbb{R}^n$. We thus try to use the characteristic method, presenting the solution by the trajectory path of the fluid particles. This is another typical way to handle transport equations. In fact, the proof in this section also works for $p=2$ as well as the case of $\Omega=\mathbb{R}^n$. Thus it is an improvement of the former section. We state a priori estimates in the case of dimensions two and three separately.
\begin{Proposition}\label{psp1}
Let $\Omega$ be a bounded smooth domain in $\mathbb{R}^2$ and $p\in (1,\infty)$, $s>\frac{2}{p} + 1$ are real numbers. Suppose $(\bs{u},\pi)$ is a smooth solution to the Euler equations with smooth initial data $\bs{u}_0\in W^{s,p}(\Omega)$. Then for any fixed $T<\infty$, there exists a constant depending only on $s,p$ and $T,\Omega$ such that for any $t\in [0,T]$,
\be\label{esp1}
\|\bs{u}(t,\cdot)\|_{s,p} + \|\nabla\pi(t,\cdot)\|_{s-1,p} \leq C \|\bs{u}_0\|^2_{s,p}.
\ee
\end{Proposition}

\begin{Proposition}\label{psp2}
Let $\Omega$ be a bounded smooth domain in $\mathbb{R}^3$ and $p\in (1,\infty)$, $s>\frac{3}{p} + 1$ are real numbers. Suppose $(\bs{u},\pi)$ is a smooth solution to the Euler equations with smooth initial data $\bs{u}_0\in W^{s,p}(\Omega)$. Then there exist a time $T^*<\infty$ depending on $\|\bs{u}_0\|_{s,p}$ and a constant $C$ depending only on $s,p$ and $T^*,\Omega$ such that for any $t\in [0,T^*]$,
\be\label{esp2}
\|\bs{u}\|_{s,p} + \|\nabla\pi\|_{s-1,p} \leq C.
\ee
\end{Proposition}

The rest of this section is devoted to the proof of these two propositions. In the following sections, we always assume that $(\bs{u},\pi)$ is the smooth solution to the initial-boundary value problem of Euler equations in a time interval $[0,T]$.
\subsection{Estimates for the Velocity and Pressure}
We first give estimates for the velocity in terms of the vorticity. In $2D$ case we introduce the stream function $\psi$ such that
\[
\bs{u}(t,x)= \nabla^\bot\psi := (\partial_{x_2}\psi(t,x), -\partial_{x_1}\psi(t,x)).
\]
Then
\[
\begin{split}
&-\Delta\psi = \omega \quad  \text{in}\ \Omega,\\
& \psi(t,x)=0 \quad \text{on} \ \Gamma.
\end{split}
\]
It follows from the elliptic estimates that
\[
\begin{split}
   & \|\psi(t,\cdot)\|_{s+1,p} \leq C \|\omega(t,\cdot)\|_{s-1,p}, \\ 
   & \|\nabla^2\psi(t,\cdot)\|_{*}\leq C\|\omega(t,\cdot)\|_{*}\leq C\|\omega(t,\cdot)\|_{\infty}\leq C
\end{split}
\]
with some constant $C$ independent of $t$ and $\omega$. Consequently,
\be\label{esp01}
\begin{split}
&\|\bs{u}(t,\cdot)\|_{s,p} \leq C \|\omega(t,\cdot)\|_{s-1,p}, \\
&\|\nabla\bs{u}(t,\cdot)\|_{*}\leq C.
\end{split}
\ee
In $3D$ case we note that
\[
\mathrm{curl}\bs{u} = \bs{\omega}, \quad \mathrm{div}\bs{u}=0 \ \text{in} \ \Omega, \quad \bs{u}\cdot\bs{n}|_{\Gamma}=0.
\]
By elliptic estimates again we have
\be\label{esp02}
\|\bs{u}(t,\cdot)\|_{s,p} \leq C \|\omega(t,\cdot)\|_{s-1,p}.
\ee
Since the pressure $\pi$ satisfies
\[
\begin{split}
&-\Delta\pi = \dv(\bs{u}\cdot\nabla\bs{u}) \text{ in }\Omega, \\
&\frac{\partial\pi}{\partial\bs{n}}|_{\Gamma} = \bs{u}\cdot\nabla\bs{u}\cdot\bs{n}.
\end{split}
\]
We have
\be\label{esp03}
\|\pi(t,\cdot)\|_{s,p}\leq C\|\bs{u}(t,\cdot)\|^2_{s,p}.
\ee

\subsection{Estimates for the Trajectory}
The trajectory ${X}(t,t_0,x)$ of the fluid particles starting from $x$ at time $t_0$ is defined as
\be\label{esp3}
\frac{{\rm d}{X}(t,t_0,x)}{{\rm d}t} = \bs{u}(t,t_0,X), \, {X}(t_0,t_0,x)=x.
\ee
In the integral form,
\be\label{esp4}
{X}(t,t_0,x)= x + \int_{t_0}^t \bs{u}(\tau, X(\tau,t_0,x)){\rm d}\tau.
\ee
Denote $X(t,x)= X(t,0,x)$. Taking derivatives with respect to $x$ in the above equation for $t_0=0$ we have
\[
\nabla_x {X}(t,x) = \mathbb{I} + \int_0^t \nabla\bs{u}(\tau, {X}(\tau,x))\nabla_x {X}(\tau,x) {\rm d}\tau.
\]
Here $\mathbb{I}$ is the identity map in $\mathbb{R}^n,\, n=2,3$. It follows that
\[
\|\nabla{X}(t,\cdot)\|_{\infty} \leq 1 + \int_0^t \|\nabla\bs{u}(\tau,\cdot)\|_{\infty}\|\nabla{X}(\tau,\cdot)\|_{\infty}{\rm d}\tau,
\]
which yields 
\be\label{esp00}
\|\nabla{X}(t,\cdot)\|_{\infty} \leq \exp(\int_0^t \|\nabla{\bs{u}}(\tau,\cdot)\|_{\infty}{\rm d}\tau),
\ee
and we also have
\be\label{esp51}
\begin{split}
&\|\nabla_x{X}(t,\cdot)\|_{s-1,p}\\ 
\leq & 1 + \int_0^t \|(\nabla\bs{u})(\tau,X(\tau,\cdot))\|_{s-1,p} \|\nabla_x {X}(\tau,\cdot)\|_{s-1,p} {\rm d}\tau \\
\leq & 1 + \int_0^t \|(\nabla\bs{u})(\tau,\cdot)\|_{s-1,p} \|\nabla_x {X}(\tau,\cdot)\|^2_{s-1,p} {\rm d}\tau\\
\leq & 1 + C(|\Omega|) \int_0^t \|\bs{\omega}(\tau,\cdot)\|_{s-1,p} \|\nabla_x {X}(\tau,\cdot)\|^2_{s-1,p} {\rm d}\tau,
\end{split}
\ee
according to the multiplication and composition properties of $W^{s,p}(\Omega)$ for $s>\frac{n}{p}$ and the elliptic estimate (\ref{esp01}) or (\ref{esp02}). Corresponding estimates hold true for $X^{-1}(t,x)={X}(0,t,x)$.

\subsection{Estimates for the Vorticity}
We first consider the evolution of the vorticity in $2D$ case, which needs some special treatments. By our assumption,
\[
s>\frac{2}{p}, \ 1<p<\infty.
\]
Hence $s>2$. We assume
\[
2<s<3.
\]
Once the estimate (\ref{esp1}) is proved in this case, it will be obvious for $s>3$. In the two dimension space, the vorticity $\omega$ is a scalar function and satisfies
\[
\begin{split}
    &\partial_t \omega + \bs{u}\cdot\nabla\omega = 0, \\
    &\omega(x,0) = \omega_0(x) :=\partial_1 {u}_{02}-\partial_2 {u}_{01}.
\end{split}
\]
By using the trajectory $X(t,x)$, we can represent $\omega$ as
\[
\omega(t,x)=\omega_0(X^{-1}(t,x)).
\]
It follows that
\be\label{esp05}
\|\omega\|_{\infty} \leq \|\omega_0\|_\infty \leq C\|\bs{u}_0\|_{s,p}.
\ee
On the other hand, by applying Proposition \ref{prop29} we have
\be\label{esp22}
\|\omega(t,\cdot)\|_{s-1,p} \leq \|\omega_0\|_{s-1,p} \exp \left( C\int_0^t \|\nabla\bs{u}(\tau,\cdot)\|_{\infty} {\rm d}\tau \right).
\ee
While by Proposition \ref{prop235} and (\ref{esp01}), (\ref{esp05}),
\[
\begin{split}
\|\nabla\bs{u}(\tau,\cdot)\|_{\infty}&\leq C\|\omega(\tau,\cdot)\|_{\infty} \ln(1+\|\omega(\tau,\cdot)\|_{s-1,p})\\
&\leq C\|\omega_0\|_{\infty} \ln(1+\|\omega(\tau,\cdot)\|_{s-1,p}).
\end{split}
\]
Substituting this into (\ref{esp22}) gives rise to
\[
\|\omega(t,\cdot)\|_{s-1,p} \leq \|\omega_0\|_{s-1,p} \exp \left( C\|\omega_0\|_{\infty}\int_0^t \ln(1+\|\omega(\tau,\cdot)\|_{s-1,p}) {\rm d}\tau \right).
\]
It follows by Gr\"{o}nwall's argument that
\be\label{esp71}
\|\omega(t,\cdot)\|_{s-1,p} \leq C\|\omega_0\|_{s-1,p}\exp\left(\exp(C\|\omega_0\|_{s-1,p}T)\right), \ \text{for any} \ t\in [0,T].
\ee
Consequently,
\be\label{esp72}
\|\bs{u}(t,\cdot)\|_{s,p}\leq C\|\bs{u}_0\|_{s,p}\exp\left(\exp(C\|\bs{u}_0\|_{s,p}T)\right), \ \text{for any} \ t\in [0,T].
\ee
By elliptic estimates (\ref{esp03}),
\be\label{esp73}
\|\pi(t,\cdot)\|_{s,p}\leq C(T)\|\bs{u}_0\|^2_{s,p}, \ \text{for any} \ t\in [0,T].
\ee
Thus the proof of Proposition \ref{psp1} is concluded.

In $3D$ case the evolution of the vorticity is
\[
\partial_t\bs\omega + \bs{u}\cdot\nabla\bs{\omega} = \bs{\omega}\cdot\nabla\bs{u}, \, \bs\omega(x,0) = \bs{\omega}_0 := \nabla \times \bs{u}_0.
\]
Along the trajectory,
\[
\partial_t\bs{\omega}|_{(t,X(t,x))} = \bs{\omega}\cdot\nabla\bs{u}|_{(t,X(t,x))}, \, \bs\omega(x,0) = \bs{\omega}_0.
\]
Hence,
\[
\bs\omega(t,X(t,x)) = \nabla{X}(t,x)\bs\omega_0(x).
\]
By the multiplication and composition property, we find
\be\label{esp9}
\begin{split}
\|\bs\omega(t,\cdot)\|_{s-1,p} & \leq C\|\bs\omega(t,X(t,\cdot))\|_{s-1,p}\|\nabla{X}^{-1}(t,\cdot)\|_{s-1,p}\\
&\leq  C\|\bs\omega_0\|_{s-1,p}\|\nabla{X}(t,\cdot)\|_{s-1,p}\|\nabla{X}^{-1}(t,\cdot)\|_{s-1,p}.
\end{split}
\ee
Using (\ref{esp51}) we find for
\[
g(t)= 1 + \|X^{-1}(t,\cdot)\|_{s-1,p} +\|X(t,\cdot)\|_{s-1,p},
\]
we have
\be\label{esp91}
g(t) \leq 1 + C \int_0^t g^4(\tau){\rm d}\tau.
\ee
Gr\"{o}nwall type argument shows that there exists a time $T^*>0$ depending on $\|\bs\omega_0\|_{s-1,p}$ and a constant $C(T^*)=C(\|\bs\omega_0\|_{s-1,p})$ such that
\[
g(t)\leq C(T^*).
\]
The estimate for $\bs{\omega}$ follows from (\ref{esp9}) and Proposition \ref{psp2} follows as before by elliptic estimates.

\noindent{\bf Remark:} More careful calculation should result in a Gr\"{o}nwall's type inequality as
\[
g(t) \leq 1 + C \int_0^t g^2(\tau){\rm d}\tau.
\]
That is, the reasonable exponential on $g(t)$ should be $2$ provided one choose suitable $g$. We do not pursue this here since we want to give a simplified proof by using characteristic method.

\section{Existence of Solutions}

The existence proof is based on an approximation argument. In fact, due to the existence results of $W^{k,p}(\Omega)$ solution with $k$ being an integer, the main problem is to construct the approximate initial data. Hence in the following, we first focus on this problem.

For completeness here we give the classical {\bf Aubin-Lions lemma} as follows:
\begin{Lemma}[Aubin-Lions Lemma]
Suppose that $G \subset H \subset K$ where G, H, and K are reflexive Banach spaces and the embedding $G\subset H$ is compact. If the sequences $u_n$ is bounded in $L^q(0,T;G), q\geq1$, and $\partial_t u_n$ is bounded in $L^p(0,T;K)$, $p\geq1$, then there exists a subsequence of $u_n$ that is strongly convergent in $L^q(0, T; H)$.
\end{Lemma}
For the proof for this lemma please refer to \cite{Rob&Rod&Sad}.

Here we start our {\bf proof} for the main theorem.

For notation convenience let us temporarily denote 
$$
\bs{u}(x)=\bs{u}_0(x)\in W^{s,p}(\Omega)
$$ 
as the given initial data, which satisfies
\[
\dv\bs{u} = 0 \text{ in }\Omega, \quad \bs{u}\cdot\bs{n}=0\ \text{ on } \Gamma.
\]
Let $\widetilde{\bs{u}}$ be the extension of $\bs{u}(x)$ in the whole space $\mathbb{R}^n, \ n=2,3$ satisfying
\be\label{51}
\dv \widetilde{\bs{u}}=0, \quad \|\widetilde{\bs{u}}\|_{s,p} \leq C \|\bs{u}\|_{s,p}.
\ee
In fact, the existence of such an extension has been proved in \cite{KMPT}.
Then we regularize $\widetilde{\bs{u}}$ as $\widetilde{\bs{u}}_{\varepsilon}$ in the whole space $\mathbb{R}^n$ so that
\be\label{52}
\wdt{\bs{u}}_{\ve}\in C^{\infty}(\mathbb{R}^n),\quad
\lim_{\ve\rightarrow 0}\|\wdt{\bs{u}}_{\ve}-\wdt{\bs{u}}\|_{s,p} = 0.
\ee
Now let $\overline{\bs{u}}_{\ve}$ be the restriction of $\widetilde{\bs{u}}_{\ve}$ on $\overline\Omega$. Then $\overline{\bs{u}_{\ve}} \in C^\infty(\overline{\Omega})$ and
\be\label{53}
\lim_{\ve\rightarrow 0} \|\overline{\bs{u}}_{\ve} - \bs{u} \|_{s,p} = 0
\ee
according to (\ref{52}). The only problem exists on the boundary value of $\overline{\bs{u}}_\ve$, which do not necessarily satisfy
\[
\overline{\bs{u}}_{\ve} \cdot \bs{n} = 0 \quad \text{on} \ \Gamma.
\]
To this end let us first note that $\overline{\bs{u}}_{\ve}$ do satisfy
\be\label{54}
\int_\Gamma \overline{\bs{u}}_{\ve}\cdot\bs{n} {\rm d}\sigma (x) = 0
\ee
due to the divergence-free condition. Now consider the Neumann problems
\[
\Delta q_{\ve} = 0 \quad \text{in} \ \Omega, \qquad \frac{\partial q_\ve}{\partial\bs{n}} = \overline{\bs{u}}_{\ve}\cdot\bs{n} \quad \text{on} \ \Gamma,
\]
which are solvable due to the compatible condition (\ref{54}). Moreover, since
\[
\bs{u}\cdot\bs{n} = 0 \quad\text{ on } \Gamma,
\]
we have
\[
\frac{\partial q_\ve}{\partial\bs{n}} = \left(\overline{\bs{u}}_{\ve}-\bs{u}\right)\cdot\bs{n} \quad \text{on} \ \Gamma.
\]
It follows by classical estimates for Neumann problems that
\be\label{55}
\|\nabla q_{\ve}\|_{s,p} \leq C \|\overline{\bs{u}}_{\ve}-\bs{u}\|_{W^{s-\frac{1}{p},p}(\Gamma)} \leq \|\overline{\bs{u}}_{\ve}-\bs{u}\|_{s,p} \rightarrow 0, 
\ee
when $\ve\rightarrow 0$.

Now we define
\[
\bs{u}_\ve = \overline{\bs{u}}_{\ve} - \nabla{q}_{\ve}.
\]
It follows that $\bs{u}_{\ve} \in C^{\infty}(\overline{\Omega})$ and
\[
\bs{u}_\ve \cdot \bs{n} = 0\quad \ \text{on} \ \Gamma.
\]
Moreover, according to (\ref{53}) and (\ref{55}), there holds
\be\label{56}
\|\bs{u}_{\ve} - \bs{u} \|_{s,p} \leq C \left( \|\overline{\bs{u}}_{\ve} - \bs{u} \|_{s,p} + \|\nabla{q}_{\ve}\|_{s,p} \right) \rightarrow 0
\ee
when $\ve\rightarrow 0$.

This completes the construction of approximate initial data.

Now we briefly give the existence of the solution to the initial boundary value problem (\ref{1.1})-(\ref{1.3}). For given $\bs{u}_0\in W^{s,p}(\Omega)$ we can find, according to the last step, $\bs{u}_{0n}\in C^{\infty}(\overline{\Omega}),\ n=1,2,3,\cdots$ such that
\[
\lim_{n \rightarrow 0} \|\bs{u}_{0n} - \bs{u}_0 \|_{s,p} = 0.
\]
Let $k=[s]+1$ and consider (\ref{1.1})-(\ref{1.3}) with the initial data $\bs{u}_0$ replaced by $\bs{u}_{0n}$. By the result of Temam(\cite{Tem}) or Bourguignon-Brezis(\cite{Bou&Bre}), there exists $T^n_*>0$ depending on $\|\bs{u}_{0n}\|_{k,p}$ and a unique solution
\[
\bs{u}_n, \ \pi_n \in C([0,T_*^n),W^{k,p}(\Omega)).
\]
According to the a priori estimates in Section 4,
\be\label{58}
\|\bs{u}_n(t,\cdot),\pi_n(t,\cdot)\|_{s,p} \leq C(\|\bs{u}_{0n}\|_{s,p}) \leq C(\|\bs{u}_0\|_{s,p}),
\ee
for any
\[
t\in (0,T^*]
\]
with $T^*$ independent of $n$. If $T^n_* < T^*$, which is mostly the case, we can extend $(\bs{u}_n(t),\pi_n(t))$ to the time interval $[0,T^*]$, by finite steps for every $n$, such that
\[
\bs{u}_n,\ \pi_n \in C([0,T^*],W^{s,p}(\Omega)), \quad \|\bs{u}_n(t,\cdot),\pi_n(t,\cdot)\|_{s,p} \leq C(\|\bs{u}_0\|_{s,p}).
\]
By taking limit in $n$ we obtain some $(\bs{u}(t,x),\pi(t,x))$ as the solution to (\ref{1.1})-(\ref{1.3}), at least satisfying, for any $\tau<s$
\[
\bs{u}(t,x),\pi(t,x) \in L^{\infty}\big(0,T^*,W^{s,p}(\Omega)\big) \cap C\big([0,T^*],W^{\tau,p}(\Omega)\big),
\]
which is due to the Aubin-Lions Lemma.

In fact, a further argument can show that
\[
\bs{u}(t,x),\pi(t,x) \in  C\big([0,T^*],W^{s,p}(\Omega)\big)
\]
and we omit the details here, please refer to \cite{Maj&Ber}. Finally, in $2D$ case, one can extend the solution in the whole time interval $[0,\infty)$ by Proposition \ref{psp1}. In $3D$ case, we can find a critical time value, still denoted as $T^*$, such that
\[
\bs{u}(t,x),\pi(t,x) \in C\big([0,T^*),W^{s,p}(\Omega)\big)
\]
by Proposition \ref{psp2}. We thus conclude the proof of our main results.

\medskip

\medskip

\end{document}